\documentclass{article}    
\usepackage{amssymb, latexsym, amsmath, amsthm}

\newcommand{\bi} {\begin{itemize}}
\newcommand{\ei} {\end{itemize}}
\newcommand{\bea} {\begin{eqnarray}}
\newcommand{\eea} {\end{eqnarray}}
\newcommand{\be} {\begin{equation}}
\newcommand{\ee} {\end{equation}}
\newcommand{\bean} {\begin{eqnarray*}}
\newcommand{\eean} {\end{eqnarray*}}

\RequirePackage{epsfig}
\RequirePackage{amssymb}
\RequirePackage{natbib}
\RequirePackage{graphics}
\bibpunct{(}{)}{,}{a}{,}{,}

\begin{document}
\title{Spectral Picture For Rationally Multicyclic Subnormal Operators}
\author{Liming Yang}
\author{
 Liming Yang \\
 Department of Mathematics \\
 Virginia Polytechnic and State University \\
 Blacksburg, VA 24061 \\
yliming@vt.edu
\date{}
     }         

\maketitle

\newtheorem{Theorem}{Theorem}
\newtheorem*{MTheorem}{Main Theorem}
\newtheorem*{TTheorem}{Thomson's Theorem}
\newtheorem*{Conjecture}{Conjecture}
\newtheorem{Corollary}{Corollary}
\newtheorem{Definition}{Definition}
\newtheorem{Assumption}{Assumption}
\newtheorem{Lemma}{Lemma}
\newtheorem{Problem}{Problem}
\newtheorem{Example}{Example}
\newtheorem*{Remark}{Remark}
\newtheorem{KnownResult}{Known Result}
\newtheorem{Algorithm}{Algorithm}
\newtheorem{Property}{Property}
\newtheorem{Proposition}{Proposition}


\abstract{
For a pure bounded rationally cyclic subnormal operator  $S$ on a separable complex Hilbert space $\mathcal H,$ \cite{ce93} shows that $clos(\sigma (S) \setminus \sigma_e (S)) = clos(Int (\sigma (S))).$
This paper examines the property for rationally multicyclic (N-cyclic) subnormal operators. We show: 
(1) There exists a 2-cyclic irreducible subnormal operator $S$ with $clos(\sigma (S) \setminus \sigma_e (S)) \neq clos(Int (\sigma (S))).$ 
(2) For a pure rationally $N-$cyclic subnormal operator $S$ on $\mathcal H$ with the minimal normal extension $M$ on $\mathcal K \supset \mathcal H,$ let $\mathcal K_m = clos (span\{(M^*)^kx: ~x\in\mathcal H,~0\le k \le m\}.$ Suppose $M | _{\mathcal K_{N-1}}$ is pure, then $clos(\sigma (S) \setminus \sigma_e (S)) = clos(Int (\sigma (S))).$       
}

\section{Introduction}

Let $\mathcal H$ be a separable complex Hilbert space and let $\mathcal L(\mathcal H )$ be the space of bounded linear operators on $\mathcal H.$ An operator $S\in \mathcal L (\mathcal H)$ is subnormal if there exist a separable complex Hilbert space $\mathcal K$ containing $\mathcal H$ and a normal operator $M_z \in \mathcal L(\mathcal K)$ such that $M_z\mathcal H \subset \mathcal H$ and $S = M_z|_{\mathcal H}.$ By the spectral theorem of normal operators, we assume that 
\[
\ \mathcal K = \oplus _{i = 1}^m L^2 (\mu _i) \tag{1-1}
\]
where $\mu_1 >> \mu_2 >> ... >>\mu_m$ ($m$ may be $\infty$) are compactly supported finite positive measures on the complex plane
$\mathbb {C},$ and $M_z$ is multiplication by $z$ on $\mathcal K.$ 
For $H= (h_1,...,h_m)\in \mathcal K$ and $G= (g_1,...,g_m)\in \mathcal K,$ we define
 \[
 \ \left\langle H(z), G(z) \right\rangle = \sum _{i = 1}^m h_i(z) \overline{g_i(z)} \dfrac{d\mu_i}{d\mu _1}, ~ | H(z) |^2 = \left\langle H(z), H(z) \right\rangle . \tag{1-2}
 \]
The inner product of $H$ and $G$ in $\mathcal K$ is defined by
 \[
 \ (H,G) = \int \left\langle H(z), G(z) \right\rangle d\mu_1(z).\tag{1-3}
 \]
$M_z$ is the minimal normal extension if 
\[
 \ \mathcal K = clos\left( span (M_z^{*k}x: ~x\in \mathcal H,~ k\ge 0)\right ). \tag{1-4}
 \]
We will always assume that $M_z$ is the minimal normal extension of $S$ and $\mathcal K$ satisfies (1-1) to (1-4). For details about the functional model above and basic knowledge of subnormal operators, the reader shall consult Chapter II of the book \cite{conway}.

For $T \in \mathcal L(\mathcal H )$, we denote by $\sigma (T)$ the
spectrum of $T,$ $\sigma _e (T)$ the essential spectrum
of $T,$ $T^*$ its adjoint, $ker(T)$ its kernel, and $Ran(T)$ its range. For a
subset $A \subset \mathbb C,$ we set $Int(A)$ for its interior, $clos (A)$  for its closure, $A^c$ for its complement, and $\bar A = \{\bar z:~ z\in A\}$. For $\lambda\in \mathbb C$ and $\delta > 0,$ we set $B(\lambda, \delta) = \{z: |z - \lambda  | <\delta \}$ and $\mathbb D = B(0,1).$ Let $\mathcal{P}$ denote the set of polynomials in the complex variable $z.$ For a compact subset $K\subset \mathbb C,$ let $Rat(K)$ be the set of all rational functions with poles off $K$ and let $R(K)$ be the uniform closure of $Rat(K).$

A subnormal operator $S$ on $\mathcal H$ is pure if for every non-zero invariant subspace $I$ of $S$ ($SI \subset I$), the operator $S|_I$ is not normal. 
For $F_1, F_2,..., F_N\in \mathcal H,$ let
 \[
 \ R^2(S | F_1,F_2,...,F_N) = clos\{r_1(S)F_1 + r_2(S)F_2 + ... + r_N(S)F_N\}
 \]
 in $\mathcal H,$ where $r_1,r_2,...,r_N \in Rat(\sigma(S))$ and let
 \[
 \ P^2(S | F_1,F_2,...,F_N) = clos\{p_1(S)F_1 + p_2(S)F_2 + ... + p_N(S)F_N\}
 \] 
in $\mathcal H,$ where $p_1,p_2,...,p_N \in \mathcal P.$
A subnormal operator $S$ on $\mathcal H$ is rationally multicyclic ($N-$cyclic) if there are $N$ vectors $F_1,F_2,...,F_N\in \mathcal H$ such that 
 \[
 \ \mathcal H = R^2(S | F_1,F_2,...,F_N)
 \]
and for any $G_1,...,G_{N-1}\in \mathcal H,$
\[
 \ \mathcal H \ne R^2(S | G_1,G_2,...,G_{N-1}).
 \]
We call $N$ is the rationally cyclic multiple of $S.$
$S$ is multicyclic ($N-$cyclic) if 
 \[
 \ \mathcal H = P^2(S | F_1,F_2,...,F_N)
 \]
and for any $G_1,...,G_{N-1}\in \mathcal H,$
\[
 \ \mathcal H \ne P^2(S | G_1,G_2,...,G_{N-1}).
 \]
We call $N$ is the cyclic multiple of $S.$
In this case, $m \le N$ where $m$ is as in (1-1). 

Let $\mu$ be a compactly supported finite positive measure on the complex plane
$\mathbb {C}$ and let $spt(\mu)$ denote the support of $\mu .$ For a compact subset $K$ with $spt(\mu) \subset K,$ let $R^2(K, \mu)$ be the closure of $Rat(K)$ in $L^2(\mu).$
 Let $P^2(\mu)$ denote the closure of
$\mathcal{P}$ in $L^2(\mu).$ 

If $S$ is rationally cyclic, then $S$ is unitarily equivalent to multiplication by $z$ on $R^2(\sigma (S), \mu _1),$ where $m=1$ and $F_1 = 1.$ We may write $R^2(S | F_1) = R^2(\sigma (S), \mu _1).$ If $S$ is cyclic, then $S$ is unitarily equivalent to multiplication by $z$ on $P^2(\mu _1).$  We may write $P^2(S | F_1) = P^2(\mu _1).$

For a rationally $N-$cyclic subnormal operator $S$ with cyclic vectors $F_1,F_2,...,F_N$ and $\lambda \in \sigma (S),$  we denote the map
 \[
 \ E(\lambda ) : \sum_{i=1}^N r_i(S)F_i \rightarrow \left[ \begin{array}{c}
r_1(\lambda )\\
r_2(\lambda ) \\
... \\
r_N(\lambda ) \end{array} \right],\tag{1-5}
 \]
where $r_1,r_2,...,r_N\in Rat(\sigma(S)).$ If $E(\lambda)$ is bounded from $\mathcal K$ to $(\mathbb C^N, \|.\|_{1,N}),$ where $\|x\|_{1,N} = \sum_{i=1}^N |x_i|$ for $x\in \mathbb C^N,$
 then every component in the right hand side extends to a bounded linear functional on $\mathcal H$ and we will call $\lambda$ a
bounded point evaluation for $S.$ We use $bpe(S)$ to denote the set of bounded point evaluations for $S.$  The set $bpe(S)$ does not depend on the choices of cyclic vectors $F_1,F_2,...,F_N$ (see Corollary 1.1 in \cite{moz}).
A point $\lambda_0 \in int(bpe(S))$ is called an analytic bounded point evaluation for $S$ if there is a neighborhood $B(\lambda_0, \delta)\subset bpe(S)$ of $\lambda_0$ such that $E(\lambda)$ is analytic as a function of $\lambda$ on $B(\lambda_0, \delta)$ (equivalently (1-5) is uniformly bounded for  $\lambda\in B(\lambda_0, \delta)$).  We use $abpe(S)$ to denote the set of analytic bounded point evaluations for $S.$  The set $abpe(S)$ does not depend on the choices of cyclic vectors $F_1,F_2,...,F_N$ (also see Remark 3.1 in \cite{moz}). Similarly, for an $N-$cyclic subnormal operator $S,$ we can define $bpe(S)$ and $abpe(S)$ if we replace $r_1,r_2,...,r_N\in Rat(\sigma(S))$ in (1-5) by $p_1,p_2,...,p_N\in \mathcal P.$

For $N = 1,$ \cite{thomson} proves a remarkable structural theorem for $P^2(\mu ).$ 
\newline

\begin{TTheorem}
There is a Borel partition $\{\Delta_i\}_{i=0}^\infty$ of $spt\mu$ such that the space $P^2(\mu |_{\Delta_i})$ contains no nontrivial characteristic functions and
 \[
 \ P^2(\mu ) = L^2(\mu |_{\Delta_0})\oplus \left \{ \oplus _{i = 1}^\infty P^2(\mu |_{\Delta_i}) \right \}.
 \]
Furthermore, if $U_i$ is the open set of analytic bounded point evaluations for
$P^2(\mu |_{\Delta_i})$ for $i \ge 1,$ then $U_i$ is a simply connected region and the closure of $U_i$ contains $\Delta_i.$
\end{TTheorem}

\cite{ce93} extends some results of Thomson's Theorem to the space $R^2(K,\mu ),$ while \cite{b08} expresses $R^2(K,\mu )$ as a direct sum that includes both Thomson's theorem and results of \cite{ce93}. For a compactly supported complex Borel measure $\nu$ of $\mathbb C,$ by estimating analytic capacity of the set $\{\lambda: |\mathcal C\nu (\lambda)| \ge c \},$ where $\mathcal C\nu$ is the Cauchy transform of $\nu$ (see Section 3 for definition), \cite{b06}, \cite{ars}, and \cite{ARS10} provide interesting alternative proofs of Thomson's theorem. Both their proofs rely on X. Tolsa's deep results on analytic capacity. There are other related research papers for $N=1$ in the history. For example, \cite{b79}, \cite{h79}, \cite{bm}, and \cite{yang3}, etc.

Thomson's Theorem shows in 
Theorem 4.11 of \cite{thomson} that $abpe(S) = bpe(S)$ for a cyclic subnormal operator $S$ (See also Chap VIII Theorem 4.4 in \cite{conway}).
Corollary 5.2 in \cite{ce93} proves that the result holds for rationally cyclic subnormal operators. For $N > 1,$ \cite{y17} extends the result to rationally $N-$cyclic subnormal operators.

It is shown in Theorem 2.1 of \cite{ce93} that if $S$ is a pure rationally cyclic subnormal operator, then
\[
 \ clos (\sigma (S) \setminus \sigma _e (S)) = clos (Int(\sigma (S))). \tag{1-6} 
 \]
This leads us to examine if (1-6) holds for a rationally $N-$cyclic subnormal operator. 

A Gleason part of $R(K)$ is a maximal set $\Omega$ in $\mathbb C$ such that for $x,y \in \Omega,$ if $e_x$ and $e_y$ denote the functionals evaluation at $x$ and $y$ respectively, then $\|e_x - e_y \| _{R(K)^*} < 2.$
\cite{ot80_2} shows that a compact set $K$ can be the spectrum of an irreducible
subnormal operator if and only if $R(K)$ has only one non-trivial Gleason part $\Omega$ and
$K = clos(\Omega).$ \cite{m88} and \cite{fm03} construct irreducible subnormal operators
with a prescribed spectrum, approximate point
spectrum, essential spectrum, and the (semi) Fredholm index.
Our first result is to construct a (rationally) 2-cyclic irreducible subnormal operator for a prescribed spectrum and essential spectrum. Consequently we show that (1-6) may not hold for a (rationally) $N-$cyclic irreducible subnormal operator with cyclic multiple $N>1.$

\begin{Theorem}\label{MTheorem1} Let $K$ and $K_e$ be two compact subsets of $\mathbb C$ such that $R(K)$ has only one nontrival Gleason part $\Omega,$ $K = clos(\Omega),$ and $\partial K \subset K_e\subset K.$ Then there exists a rationally 2-cyclic irreducible subnormal operator $S$ such that $\sigma (S) = K,$ $\sigma _e(S) = K_e,$ and $ind(S-\lambda ) = -1$ for $\lambda \in K\setminus K_e.$ If, in particular, $\mathbb C \setminus K$ has only one component, then $S$ can be constructed as a 2-cyclic irreducible subnormal operator.
\end{Theorem} 

Let $K = clos(\mathbb D)$ and $K_e = \partial \mathbb D \cup clos(\frac{1}{2}\mathbb D).$ We see that 
 \[
 \ clos(K\setminus K_e) = \{z:~ \frac{1}{2} \le |z| \le 1 \} \ne clos(Int(K)) = clos(\mathbb D).
 \]
From Theorem \ref{MTheorem1}, we get the following result.

\begin{Corollary}
There exists a 2-cyclic irreducible subnormal operator $S$ such that (1-6) does not hold.
\end{Corollary} 

In the second part of this paper, we will investigate certain classes of rationally $N-$cyclic subnormal operators that have the property (1-6). Let $S$ be a rationally $N-$cyclic subnormal operator on $\mathcal H = R^2(S | F_1,F_2,...,F_N).$ Let $\psi$ be a smooth function with compact support. Define
 \[
 \ \mathcal K _n^\psi = clos \left \{\psi ^m x:~ x\in \mathcal H, ~ 0\le m \le n \right \}, 
 \] 
then
 \[
 \ \mathcal H \subset \mathcal K _1^\psi \subset ... \subset \mathcal K _n^\psi \subset...\subset \mathcal K
 \]
and $M_z |_{\mathcal K _n^\psi}$ is a subnormal operator.

\begin{Definition} \label{SubnormalClass}
A subnormal operator satisfies the property $(N,\psi)$ if the following conditions are met:
\newline
(1) $S$ is a pure (rationally) $N-$cyclic subnormal operator on $\mathcal H = R^2(S | F_1,F_2,...,F_N).$ 
\newline
(2) $\psi$ a smooth function with compact support and $Area(\sigma(S) \cap \{\bar \partial  \psi = 0\}) = 0.$
Let $M_z$ on $\mathcal K$ be the minimal normal extension of $S$ satisfying (1-1) to (1-4), then $M_z |_{\mathcal K _{N-1}^\psi}$ is also a pure subnormal operator. 
\end{Definition}

\begin{Theorem}\label{MTheorem2}
Let $N > 1$ and let $S$ be a pure subnormal operator on $\mathcal H$ satisfying the property $(N,\psi ),$ then there exist bounded open subsets $U_i$ for $ 1\le i \le N$ such that 
\[
 \ \sigma _e (S) = \bigcup _{i=1}^N \partial U_i, ~\sigma (S) = \bigcup _{i=1}^N clos(U_i),
 \]
and 
 \[
 \ ind(S - \lambda ) = -i
 \]
for $\lambda \in U_i $ and $i = 1,2,...N.$ Consequently,
 \[
 \ \sigma (S) = clos (\sigma (S) \setminus \sigma _e (S)) = clos (Int(\sigma (S))). 
 \] 
\end{Theorem}

An important special case is that $\psi = \bar z.$ In section 3, we will provide several examples of subnormal operators that satisfy the property $(N,\psi).$ 
We prove Theorem \ref{MTheorem1} in section 2 and Theorem \ref{MTheorem2} in section 3.

\section{Spectral Pictures for Irreducible Rationally 2-Cyclic Subnormal Operators}

In this section, we assume that $K$ is a compact subset of $\mathbb C,$ $Int(K) \neq \emptyset,$ and $R(K)$ has only one nontrival Gleason part $\Omega$ with $K = clos(\Omega ).$ Theorem 5 and Corollary 6 in \cite{m88} constructs a  representing measure  $\nu$ of $R(K)$ at $z_0\in Int(K)$ with support on $\partial K$ such that $S_\nu$ on $R^2(K, \nu )$ is irreducible, $\sigma (S_\nu ) = K,$ $\sigma _e(S_\nu) = \partial K,$ and $ind(S_\nu - \lambda ) = -1$ for $\lambda \in Int(K) =\sigma (S_\nu) \setminus \sigma _e(S_\nu).$ From Theorem 6.2 in \cite{gamelin}, we get
 \[
 \ L^2 (\nu ) = R^2(K, \nu ) \oplus N^2 \oplus \overline {R_0^2(K, \nu )}\tag{2-1}
 \]
where $\overline {R_0^2(K, \nu )} = \{\bar r: r(z_0) = 0 \text{ and } r \in R^2(K, \nu )\}.$ The operator $M_z,$ multiplication by $z$ on $L^2 (\nu ),$ can be written as the following matrix with respect to (2-1):
 \[
 \ M_z = \left [ \begin{matrix} S_\nu, & A,& B \\ 0, & C,& D \\0,&0,&T^*_\nu\end{matrix} \right ]
 \]
where $T_\nu,$ multiplication by $\bar z$ on $\overline {R_0^2(K, \nu )},$ is an irreducible rationally cyclic subnormal operator with $\sigma (T_\nu ) = \bar K,$ $\sigma _e(T_\nu) = \partial \bar K,$ and $ind(T_\nu - \lambda ) = -1$ for $\lambda \in Int(\bar K).$ Let 
\[
 \ S = \left [ \begin{matrix} S_\nu, & A\\ 0, & C\end{matrix} \right ],
 \] 
then $S$ is the dual of $T_\nu.$ From the properties of dual subnormal operators (see, for example, Theorem 2.4 in \cite{fm03}), we see that $S$ is an irreducible subnormal operator with $\sigma (S) = K,$ $\sigma _e(S) = \partial K,$ and $ind(S - \lambda ) = -1$ for $\lambda \in Int(K).$  

The following lemma, due to \cite{cd78} on page 194, allows us to choose eigenvectors for $S^*$ in a co-analytic manner whenever the Fredholm index function for $S$ is $-1.$
\begin{Lemma} \label{CDLemma}
If $X \in L(\mathcal H)$ and $ind(X - \lambda ) = -1$ for all $\lambda \in G := \sigma(X)\setminus \sigma _e(X),$ then there exists a co-analytic function $h : G\rightarrow H$ that is not identically zero on
any component of $G$ such that $h(\lambda ) \in ker(X - \lambda )^*.$ In particular, for every $x \in \mathcal H,$
the function $\lambda \rightarrow (x, h(\lambda ))$ is analytic on $G.$
\end{Lemma} 

Using Lemma \ref{CDLemma}, we conclude that there exists a co-analytic function $k_\lambda \in \mathcal H := R^2(K, \nu ) \oplus N^2$ such that $(S -\lambda )^* k_\lambda = 0$ on $Int(K).$ 
Let $\delta _{\lambda }$ be the point mass measure at $\lambda.$ Let $K_e \subset K$ be a compact subset of $\mathbb C$ such that $\partial K \subset K_e.$ Let $\{\lambda _n \} \subset K_e\cap Int(K)$ with $K_e\cap Int(K) \subset clos (\{\lambda _n \}).$ 
Define
 \[
 \ \mu = \nu + \sum_{n=1}^\infty c_n \delta _{\lambda _n},\tag{2-2}
 \] 
where $c_n > 0$ and $\sum_{n=1}^\infty c_n \|k_{\lambda _n} \|^2 = 1.$ Let $M_z^1$ be the multiplication by $z$ operator on $L^2(\mu ).$

\begin{Lemma}\label{DualLemma} Define an operator $T$ from $\mathcal H$ to $L^2(\mu )$ by
 \[
 \ Tf(z) = \begin{cases} f(z),& z\in \partial K\\ (f,k_{\lambda _n}),&z=\lambda _n. \end{cases} \tag{2-3}
 \]
Then $T$ is a bounded linear one to one operator with closed range. Set $\mathcal H_1 = Ran(T),$ then $T$ is invertible from $\mathcal H$ to $\mathcal H_1,$ $M_z^1\mathcal H_1\subset \mathcal H_1,$ $S_1 = M_z^1 | _{\mathcal H_1}$ is an irreducible subnormal operator such that $S_1= TST^{-1},$ and $M_z^1$ is the minimal normal extension of $S_1.$
\end{Lemma} 

{\bf Proof:} By definition, we get
 \[
 \ \|f\|_{L^2(\nu )}^2 \le \|Tf\|_{L^2(\mu )}^2 = \|f\|_{L^2(\nu )}^2 + \sum_{n=1}^\infty c_n |(f,k_{\lambda _n})|^2 \le 2\|f\|_{L^2(\nu )}^2.
 \]
Therefore, $T$ is a bounded linear operator and invertible from $\mathcal H$ to $\mathcal H_1.$ Since $(zf,k_{\lambda _n}) = \lambda _n (f,k_{\lambda _n}),$ we see that $M_z^1\mathcal H_1\subset \mathcal H_1$ and $S_1= TST^{-1}.$ Since $(Tk_{\lambda _n}) (\lambda _n) = \|k_{\lambda _n}\|^2 > 0,$ clearly, we have 
 \[
 \ L^2(\mu ) = clos\left (span\{\bar z^mx: ~x\in\mathcal H_1,~ m \ge 0\}\right ).
 \]
Therefore, $M_z^1$ is the minimal normal extension of $S_1.$

It remains to prove that $S_1$ is irreducible. Let $N_1$ and $N_2$ be two reducing subspaces of $S_1$ such that $\mathcal H_1 = N_1\oplus N_2.$ Then for $f_1\in N_1$ and $f_2\in N_2,$ we have
 \[
 \ (z^n f_1, z^mf_2) = \int z^n\bar z^m f_1\bar f_2 d\mu = 0
 \]
for $n,m = 0,1,2,....$ This implies $f_1(z)\bar f_2(z) = 0$ a.e. $\mu.$ By the definition of $T,$ we see that $(T^{-1}f_1)(z)\overline{(T^{-1}f_2)}(z) = 0$ a.e. $\nu.$ Hence, $\mathcal H = T^{-1}N_1\oplus T^{-1}N_2.$ $T^{-1}N_1$ and $T^{-1}N_2$ are reducing subspaces of $S.$ By the construction, $T_\nu$ is irreducible (Corollary 6 in \cite{m88}), so $S,$ as the dual $T_\nu,$ is irreducible (see, for example, Theorem 2.4 in \cite{fm03}). This means that $N_1 = 0$ or $N_2 = 0.$ The lemma is proved.  

We write the operator $M_z^1$ as the following:
 \[
 \ M_z^1 = \left [ \begin{matrix}S_1,& A_1\\0,& T_1^*\end{matrix} \right ] \tag{2-4} 
 \]
Then $T_1,$ as a dual of $S_1,$ is irreducible. 

\begin{Lemma} \label{GeneratorLemma}
Let $\mu$ be as in(2-2) and let $\mathcal H_1$ be as in Lemma \ref{DualLemma}. Define
 \[
 \ F(z) = \begin{cases} \bar z - \bar z_0, ~ &z\in \partial K, \\ 0, ~ &z\in Int(K).\end{cases}\tag{2-5}
 \]
and 
 \[
 \ G_n(z) = \begin{cases} k_{\lambda _n} (z), ~ &z\in \partial K, \\ -1/c_n, ~ &z= \lambda _n, \\ 0, ~ &z=\lambda _m, ~m \neq n.\end{cases}\tag{2-6}
 \]
Then
 \[
 \ \mathcal H_1 ^\perp = clos\left ( span\{r(\bar z)F,G_j, ~ 1\le j < \infty,~ r\in Rat(K) \}\right ).
 \]
\end{Lemma}

{\bf Proof: } It is straightforward to check, from (2-1), (2-2), and (2-3), that $F, G_j \in \mathcal H_1 ^\perp.$ Now let $H(z) \perp clos\left ( span\{r(\bar z)F,G_j, ~ 1\le j < \infty,~ r\in Rat(K) \}\right ),$ then
 \[
 \ \int H(z) r(z) \bar F(z) d\mu = \int H(z) r(z)(z-z_0) d\nu = 0
  \] 
	for $r\in Rat(K).$ From (2-1), we see that the function $H |_{\partial K} \in \mathcal H.$  It follows from $\int H(z) \bar G_j(z) d\mu = 0$ that $ H(\lambda _j) = (H |_{\partial K}, k_{\lambda _j}).$ Thus, $H(z )\in \mathcal H _1.$ The lemma is proved.

\begin{Lemma} \label{2CyclicLemma}
Let $\mu,$ $T_1,$ $F,$ and $G_n$ be as in (2-2), (2-4), (2-5) and (2-6), respectively. Then there exists a sequence of positive numbers $\{a_n\}$ satisfying
 \[
 \ \sum_{n = 1}^\infty a_n \|G_n\| < \infty, ~ G = \sum_{n = 1}^\infty a_nG_n, 
 \]
and
 \[
 \ \mathcal H_1^\perp = clos\left ( span \{r(\bar z ) F(z) + p(\bar z) G(z): ~r\in Rat(K), ~ p\in \mathcal P \} \right ).
 \]
Therefore, $T_1$ is a rationally 2-cyclic irreducible subnormal operator with
 \[
 \ \sigma (T_1) = \bar K,~ \sigma _e (T_1) = \bar K_e,  \text{ and } ind (T_1 - \lambda ) = -1, ~\lambda \in \bar K \setminus \bar K_e.\tag{2-7}
 \]

\end{Lemma}

{\bf Proof: } Notice that
 \[
 \ \int f(z) (z - \lambda _n) \bar k_{\lambda _n}(z) d \nu = 0
 \]
for $f\in \mathcal H.$ We conclude, from (2-1),  that $(\bar z - \bar \lambda _n) k_{\lambda _n}(z)\in \overline{R_0^2 (K, \nu )}.$ Hence, there are $\{r_n\}\subset R^2 (K, \nu )$ such that
 \[
 \ k_{\lambda _n}(z) = \dfrac{r_n(\bar z)}{\bar z - \bar \lambda _n}(\bar z - \bar z_0).
 \]
 We will recursively choose $\{a_n\}.$ First choose $a_1 = 1.$ Then we assume that $a_1,a_2,...,a_n$ have been chosen. Now we will choose $a_{n+1}.$ Let
 \[
 \ p_k(z) = \dfrac{\Pi _{j\neq k, 1\le j \le n}(z - \bar \lambda _j)}{a_k\Pi _{j\neq k, 1\le j \le n}(\bar \lambda _k - \bar \lambda _j)},
 \]
for $k = 1,2,...,n.$ Denote
 \[
 \ q_{1k}(z) = p_k(z) \sum _{j\neq k, 1\le j \le n} \dfrac{a_j}{z-\bar\lambda _j} r_j(z) 
 \]
and
\[
 \ q_{2k}(z) = \dfrac{a_k (p_k(z)  - p_k(\bar\lambda _k) )}{z-\bar\lambda _k}r_k(z). 
 \] 
So $p_{k} \in \mathcal P$ and $q_{1k}, q_{2k}\in R^2(K, \nu)$ for $k = 1,2,...,n.$ Clearly,
 \[
 \ p_k(\bar z) \sum _{j=1}^n a_j G_j(z) - (q_{1k}(\bar z) + q_{2k}(\bar z))(\bar z - \bar z_0) = \dfrac{r_k(\bar z)(\bar z - \bar z_0)}{\bar z - \bar \lambda _k}, ~ z\in \partial K.
 \]
Hence,
 \[
 \ p_k(\bar z) \sum _{j=1}^n a_j G_j(z) - (q_{1k}(\bar z) + q_{2k}(\bar z))F(z) = G_k(z),~ a.e. ~\mu.
 \]
We have the following calculation:
 \[
 \ \begin{aligned}
 \ & \int \left | p_k(\bar z) \sum _{j=1}^{n+1} a_j G_j(z) - (q_{1k}(\bar z) + q_{2k}(\bar z))F(z) - G_k(z) \right |^2 d \mu \\
 \ = & \int \left | p_k(\bar z) a_{n+1} G_{n+1}(z) \right |^2 d \mu  \\
 \ \le & \left (\dfrac{a_{n+1}}{a_k} \right )^2 \dfrac{(4D^2)^{n-1}}{\Pi _{j\neq k, 1\le j \le n}|\lambda _k - \lambda _j|^2}\|G_{n+1}\|^2 
 \ \end{aligned}
 \]
where $D = \max\{|z|: ~ z \in K\}.$
Now set
 \[
 \ a_{n+1} = \min\left ( \dfrac{1}{2^{n+1}}, \min_{1\le k \le n}\dfrac{a_k\Pi _{j\neq k, 1\le j \le n}\min(1,|\lambda _k - \lambda _j|)}{4^n\max(1,D)^{n-1}}\right ) / \max(1,\|G_{n+1}\|).\tag{2-8}
 \]
So we have chosen all $\{a_n\}.$ 
From (2-8), we have the following calculation.
 \[
 \ \begin{aligned}
 \ & \left\| p_k \sum_{i = n+2}^\infty a_j G_j \right \| \\
  \ & \le \dfrac{(2D)^{n-1}}{a_k \Pi _{j\neq k, 1\le j \le n}|\lambda _k - \lambda _j| }\sum_{i = n+2}^\infty \dfrac{a_k\Pi _{j\neq k, 1\le j \le i-1} \min(1,|\lambda _k - \lambda _j|)}{4^{i-1}\max(1,D)^{i-2}} \\
	\ & \le \dfrac{1}{2^{n+2}}.
 \ \end{aligned}
 \]
Therefore,
 \[
 \ \begin{aligned}
 \ & \left \| p_k(\bar z) G - (q_{1k}(\bar z) + q_{2k}(\bar z))F - G_k(z) \right \| \\
 \ & \le \left \| p_k(\bar z) \sum _{j=1}^{n+1} a_j G_j - (q_{1k}(\bar z) + q_{2k}(\bar z))F - G_k(z) \right \| + \left \| p_k(\bar z) \sum _{j=n+2}^{\infty} a_j G_j \right \| \\
 \ & \le \dfrac{1}{2^n}. 
 \ \end{aligned}
 \]
Hence, 
 \[
 \ G_k \in clos\left ( span \{r(\bar z ) F(z) + p(\bar z) G(z): ~r\in Rat(K), ~ p\in \mathcal P \} \right ), ~ k = 1,2,....
 \]
Since $T_1$ is the dual of $S_1,$ we see that $\sigma (M_z^1) \subset \sigma _e(S_1)\cup \overline{\sigma _e(T_1)}$ (see, for example, Theorem 2.4 in \cite{fm03}), $\sigma _e(S_1) = \partial K,$ and $\sigma _e(T_1) \supset \partial \bar K.$ So (2-7) follows. 
 This completes the proof.

{\bf Proof of Theorem \ref{MTheorem1}:} It follows from Lemma \ref{2CyclicLemma}.

\section{Spectral Picture of a Class of Rationally Multicyclic Subnormal Operators}

In this section, we will prove Theorem \ref{MTheorem2}. First we provide some examples of subnormal operators that have the property $(N,\psi)$ in Definition \ref{SubnormalClass}.

\begin{Example} \label{FRSCExample}
Every pure subnormal operator $S$ on $\mathcal H$ with finite rank self-commutator has the property $(N,\psi).$ Notice that the structure of such subnormal operators has been established based on Xia's model (see \cite{x96} and \cite{y98}).
\end{Example}

{\bf Proof:} Assume that $M_z$ on $\mathcal K$ is the minimal normal extension satisfying (1-1) to (1-4). Define the self-commutator as the following
 \[
 \ D = [S^*,S] = S^*S - SS^*. 
 \]
The element $x\in ker(D)$ if and only if $\bar z x\in \mathcal H.$ This implies $S ker (D) \subset ker(D).$ Therefore,
 \[
 \ S^* Ran(D) \subset Ran(D). \tag{3-1}
 \]
Let
 \[
 \ \mathcal H_0 = clos \left ( span (S^nf: f\in Ran(D), ~ n \ge 0)\right ),
 \]
then $S|_{\mathcal H_0}$ is N-cyclic subnormal where $N = dim (Ran(D)).$ 

On the other hand,
 \[
 \ S^*S^n D = SS^*S^{n-1}D + DS^{n-1}D,
 \]
hence, we can recursively show that $S^*S^nRan(D) \subset \mathcal H_0$ since (3-1). So $S^* \mathcal H_0 \subset \mathcal H_0.$ This implies that
 \[
 \ S (\mathcal H \ominus \mathcal H _ 0) \subset \mathcal H \ominus \mathcal H _ 0
 \] 
and $S |_{\mathcal H \ominus \mathcal H _ 0}$ is normal. Since $S$ is pure, we conclude that $\mathcal H = \mathcal H_0$ and $S$ is N-cyclic. 
From (3-1), we see that there is a polynomial $p$ such that
 \[
 \ \bar p(S^* | _{Ran(D)}) = 0.
 \]
Therefore,
\[
 \ p(S): \mathcal H \rightarrow ker(D).
 \]
Hence,
 \[
 \ \| M_z^* p(S)f \| = \| M_z p(S)f \| = \| S p(S)f \| = \| S^* p(S)f \| 
 \]
for $f\in \mathcal H.$ This implies $\bar z p\mathcal H \subset \mathcal H.$ Let $\psi = \bar z p,$ then $Area\{\bar\partial \psi = 0 \} = Area\{z: ~ p(z) = 0 \} = 0,$ $\mathcal K_{N-1}^\psi = \mathcal H,$ and $S$ satisfies the property $(N,\psi )$ in Definition \ref{SubnormalClass}.

\begin{Example}\label{TCExample}
In Lemma \ref{2CyclicLemma}, if $K=clos(\mathbb D)$ and $K_e=(\partial \mathbb D)\cup (\frac{1}{2}\partial \mathbb D),$ then the operator $T_1$ is a 2-cyclic irreducible subnormal operator satisfying the property $(2,\psi )$ where $\psi = |z|^4 - \frac{5}{4} |z|^2.$ 
\end{Example}

{\bf Proof:} For $f\in \mathcal H_1,$ we get
 \[
 \ \psi f = (|z|^2 - 1)(|z|^2 - \frac{1}{4})f + \frac{1}{4}f = \frac{1}{4}f
 \]
since $spt(\mu ) \subset K_e.$ Hence, $\mathcal K_1^\psi = \mathcal H_1.$ On the other hand, 
 \[
 \ Area\{\bar\partial \psi = 0 \} \le Area\left ( \{0\}\cup \{|z| = \frac{5}{8}\}\right) = 0.
 \]
Therefore, the operator $T_1$ satisfies the property $(2, \psi ).$

In the remaining section, we assume that $N > 1$ and $S$ is a pure rationally $N-$cyclic subnormal operator on $\mathcal H = R^2(S | F_1, F_2, ...,F_N)$ and $M_z$ on $\mathcal K,$ which satisfies (1-1) to (1-4), is the minimal normal extension of $S.$ Moreover, $S$ satisfies the property $(N,\psi )$ in Definition \ref{SubnormalClass}. Let $U_k$ be the set of $\lambda\in Int(\sigma (S))$ such that $Ran(S - \lambda )$ is closed and 
 $ dim\left ( ker(S - \lambda )^*\right ) = k,$
where $k = 1,2,...,N.$ 

\begin{Lemma} \label {KBPELemma}
If $1\le k \le N,$ $\delta >0,$ $B(\lambda_0,2 \delta) \subset Int(\sigma(S)),$ $I$ is an index subset of $\{1,2,...,N\}$ with size $N-k,$  $F = \sum_{i=1}^N r_iF_i$ where $r_i\in Rat(\sigma(S)),$ and $\{a_{ls}(\lambda )\}_{1\le l \le N-k,1\le s \le k}$ are analytic on $B(\lambda_0,2 \delta)$ such that 
 \[
 \ \sup _{1 \le s\le k,\lambda\in B(\lambda_0,\delta)} | r_{j_s}(\lambda ) + \sum_{l=1}^{N-k} a_{ls} (\lambda ) r_{i_l}(\lambda )| \le M \| F\| \tag{3-2}
 \]
and 
 \[
 \ F_{i_l}(z ) = \sum_{s=1}^{k} a_{ls} (z) F_{j_s}(z), ~ a.e ~ \mu_1 | _{B(\lambda_0,\delta)},\tag{3-3}
 \]
where $i_l\in I$ and $j_s\notin I.$ Then $\lambda _0 \in \bigcup _{i=k}^N U_k.$
\end{Lemma}

{\bf Proof:} From (3-3), we get
 \[
 \ \int_{B(\lambda_0,\delta)} |F|^2 d\mu_1 = \int_{B(\lambda_0,\delta)} \left | \sum_{s=1}^k\left (r_{j_s}(z) + \sum_{l=1}^{N-k} a_{ls} (z) r_{i_l}(z)\right ) F_{j_s}(z) \right | ^2 d \mu_1 .
 \]
Using (3-2) and the maximal modulus principle,
\[
 \ \sup _{1 \le s\le k,\lambda\in B(\lambda_0,\delta)} \left | r_{j_s}(\lambda ) + \sum_{l=1}^{N-k} a_{ls} (\lambda ) r_{i_l}(\lambda ) \right | \le\dfrac{M}{\delta}\|(S-\lambda_0) F\|.
\]
Hence,
\[
 \ \int |F|^2 d\mu_1 \le \int_{B(\lambda_0,\delta)^c} |F|^2 d\mu_1 + (\sum_{j \notin I} \| F_j\|)^2 \sup _{1 \le s\le k,\lambda\in B(\lambda_0,\delta)} \left | r_{j_s}(\lambda ) + \sum_{l=1}^{N-k} a_{ls} (\lambda ) r_{i_l}(\lambda ) \right |^2.
\]
Therefore,
\[
 \ \| F\| \le M_1 \|(S - \lambda _0) F\|,
\]
where
 \[
 \ M_1^2 = \left ( 1  + (\sum_{j \notin I} \| F_j\|)^2 \right )( \dfrac{M}{\delta} )^2.
\]
So $Ran(S - \lambda _0)$ is closed. On the other hand, there are $k$ linearly independent $k^j_\lambda\in \mathcal H$ such that
 \[
 \ r_{j_s}(\lambda ) + \sum_{l=1}^{N-k} a_{ls} (\lambda ) r_{i_l}(\lambda ) = \int \left\langle F(z), k_\lambda^j(z) \right\rangle d\mu_1(z)
 \]
where $j_s \notin I$ and $\lambda\in B(\lambda_0,\delta).$ This implies
 \[
 \ dim ( Ker(S - \lambda _0)^*) \ge k.
 \]
Therefore, $\lambda _0 \in \cup_{i=k}^N U_i.$ 

Let $\nu$ be a compactly supported finite measure on $\mathbb {C}.$ The transform
\[
\ \mathcal C_\psi^i\nu (z) = \int \dfrac{(\psi ( w) - \psi (z))^i}{w - z} d\nu (w)
\]
is continuous at each point $z$ with $|\nu | (\{z\}) = 0$ and $i > 0.$ For $i = 0,$ the transformation 
 \[
 \ \mathcal C_\psi^0 (\nu ) = \mathcal C (\nu ) = \int \dfrac{1}{w - z} d\nu (w)
 \]
 is the Cauchy transform of $\nu.$ Let $M^G(z)$ be the following $N$ by $N$ matrix,
 \[
 \ M^G(z) = \left [ \mathcal C_\psi^{i-1} (\left\langle F_j, G \right\rangle \mu_1)\right ] _{N\times N}
 \]
where we assume that $G\perp \mathcal K_{N-1}^\psi$ or equivalently $G$ satisfies the following conditions
 \[
 \ \bar\psi^n G \perp \mathcal H, ~n = 0,1,2,...,N-1. \tag{3-4}
 \]
The set $W^G\subset \mathbb C$ is defined by:
\[
\ W^G = \{\lambda:~ \int \dfrac{1}{|z - \lambda|} |\langle F_i(z),G(z)\rangle | d \mu_1 ( z) < \infty, ~ 1\le i\le N\}. 
\]
Let 
 \[
 \ \Omega^G = Int(\sigma(S))\cap W^G \cap \{ \lambda : |det(M^G(\lambda ) | > 0 \}.\tag{3-5} 
 \]
Then for $\lambda \in \Omega^G,$ the matrix
 \[
 \ \left [\mathcal C (\langle F_j \psi^{i-1}, G\rangle \mu_1 )\right ] _{N\times N}\tag{3-6}
 \]
is invertible. By Construction, we see that
 \[
 \ det(M^G(z)) = 0 ~ a.e. ~ Area | _{(clos(\Omega^G))^c}.
 \]

\begin{Lemma}\label{KeyLemma0}
Using above notations, we conclude that
 \[
 \ \Omega ^G \subset abpe(S).
 \]
Hence, by Lemma \ref{KBPELemma}, we get $\Omega ^G \subset U_N.$ 
\end{Lemma}

{\bf Proof:} Using (3-4), (3-5),  and (3-6), we see that the lemma is a direct application of Theorem 2 in \cite{y17}.

Let $A = \{\lambda _n : ~ \mu _1(\{\lambda _n \}) > 0\}$ be the set of atoms for $\mu_1.$ Now let us define the matrix $M_j^G(z)$ to be a submatrix of $M^G(z)$ by eliminating the first row and $j$ column. Let $B_j^G(z)$ be the $j$ column of the matrix $M^G(z)$ by eliminating the first row. Define
\[
 \ \Omega_j^G = \left ( Int(\sigma(S))\cap A^c\cap\{z: ~ |det(M_j^G(z)) | > 0\} \right ) \setminus clos( \Omega^G ).\tag{3-7}
 \]
Notice that $M_j^G(\lambda )$ is continuous at each $\lambda\in \Omega_j^G.$ On $\Omega_j^G,$ we can define the following function vector
 \[
 \ a_j (z) = [ a_{ij} (z) ] _{(N-1)\times 1} = (M_j^G(z))^{-1} B_j^G(z).\tag{3-8}
 \] 

\begin{Lemma}\label{KeyLemma}
Let  $G,$ $\Omega^G,$ $\Omega^G_j,$  and $a_j(z)$ be as in (3-4), (3-5), (3-7), and (3-8), respectively. 
Then for $\lambda_0\in \Omega^G_j,$ there exists $\delta >0$ such that $a_j(z)$ equals an analytic function on $B(\lambda _0, \delta) \subset Int(\sigma(S)$ almost everywhere with respect to the area measure. Moreover, 
 \[
 \ \mathcal C(\langle F_j, G\rangle \mu)(z) = \sum_{k=1}^{j-1} a_{kj}(z)\mathcal  C(\langle F_k, G\rangle \mu)(z) + \sum_{k=j+1}^N a_{k-1,j}(z)\mathcal  C(\langle F_k, G\rangle \mu)(z), ~ a.e.~ Area |_{B(\lambda _0, \delta)},\tag{3-9}
 \]
and
 \[
 \ \langle F_j, G\rangle   = \sum_{k=1}^{j-1} a_{kj}(z)\langle F_k, G\rangle  + \sum_{k=j+1}^N a_{k-1,j}(z)\langle F_k, G\rangle , ~a.e. \mu |_{B(\lambda _0, \delta)}.\tag{3-10}
 \]
\end{Lemma}

Proof: Without loss of generality, we assume that $j=N.$ For $z\in Int(\sigma(S) \cap W^G \cap \Omega^G_N,$ write
 \[
 \ M^G(z) = \left[ \begin{array}{cc}
A_N^G(z) & c_N^G(z) \\
M_N^G(z) & B_N^g(z) \end{array} \right]
 \]
where
 \[
 \ A_N^G(z) = [ \mathcal C(\langle F_1, G\rangle \mu_1)(z), \mathcal C(\langle F_2, G\rangle \mu_1)(z), ... , \mathcal C(\langle F_{N-1}, G\rangle \mu_1)(z) ]
 \]
and
 \[
 \ c_N^G(z) = \mathcal C(\langle F_N, G\rangle \mu_1)(z).
 \]
By construction of $\Omega^G_N,$ we conclude that
 \[
 \ det(M^G(z)) = (A_N^G(z)(M_N^G(z))^{-1}B_N^G(z) - c_N^G(z)) det(M_N^G(z)) = 0 ~ a.e. ~ Area|_{\Omega^G_N}.
 \] 
Therefore,
 \[
 \ c_N^G(z) = A_N^G(z)(M_N^G(z))^{-1}B_N^G(z)  ~ a.e. ~ Area|_{\Omega^G_N}.\tag{3-11}
 \] 

Let $\nu_i = \langle F_i, G\rangle\mu_1$ and $H_{i,m} (z) = \frac{m^2}{\pi} \nu_i(B(z, \frac{1}{m})),$ then the functions $H_{i,m} (z)$ are bounded with compact supports. We have
 \[
 \ \mathcal C(H_{i,m} dA) (w) = \int_{|\lambda - w| \ge \frac{1}{m}} \dfrac{1}{\lambda - w} d\nu_i(\lambda ) + \int_{|\lambda - w| < \frac{1}{m}} \dfrac{m^2|\lambda - w|^2}{\lambda - w} d\nu_i(\lambda ).
 \]
Hence,
\[
\ | \mathcal C(H_{i,m} dA) (w) - \mathcal C\nu_i (w) | \le 2 \int _{|w - z | <1/m} \dfrac{1}{|w - z | } d |\nu_i| (z) ~a.e.~ Area
\]
and
\[
\ \lim_{m\rightarrow \infty}\mathcal C(H_{i.m} dA) (w) =\mathcal C\nu_i (w), ~ a.e.~ Area.
\]
Let $C_0 >0 $ be a constant such that $|\psi (z) - \psi (w)|\le C_0|z-w|.$
We estimate $\mathcal  C^1_\psi(\nu_i)$ as the following,
\[
 \ \begin{aligned}
\ & | \mathcal  C^1_\psi(H_{i,m} dA) (w) - \mathcal  C^1_\psi\nu_i (w) |  \\
\ = & \left | \dfrac{m^2}{\pi} \int \int _{|z-\lambda | < \frac{1}{m}} \dfrac{\psi (z) - \psi (w)}{z - w} dA(z) d\nu_i (\lambda ) - \mathcal  C^1_\psi\nu_i (w) \right | \\
\ \le & \left | \dfrac{m^2}{\pi} \int _{|\lambda-w|\ge \frac{1}{\sqrt{m}}} \int _{|z-\lambda | < \frac{1}{m}} \dfrac{\psi (z) - \psi (w)}{z - w} dA(z) d\nu_i (\lambda ) - \int _{|\lambda-w|\ge \frac{1}{\sqrt{m}}} \dfrac{\psi(\lambda ) - \psi ( w)}{\lambda - w} d\nu_i (\lambda )\right | \\
\ + & \left | \dfrac{m^2}{\pi}  \int _{|\lambda-w| < \frac{1}{\sqrt{m}}} \int _{|z-\lambda | < \frac{1}{m}} \dfrac{\psi (z) - \psi (w)}{z - w} dA(z) d\nu_i (\lambda ) \right | +  \left | \int _{|\lambda-w|< \frac{1}{\sqrt{m}}} \dfrac{\psi(\lambda ) - \psi ( w)}{\lambda - w} d\nu_i (\lambda ) \right | 
\ \end{aligned}
 \]
Notice that 
 \[
 \ \dfrac{m^2}{\pi}\int _{|\lambda-w| \ge \frac{1}{\sqrt{m}}} \int _{|z-\lambda | < \frac{1}{m}} \dfrac{1}{z - w} dA(z) d\nu_i (\lambda ) = \int _{|\lambda-w| \ge \frac{1}{\sqrt{m}}} \dfrac{1}{\lambda - w} d\nu_i (\lambda ).
 \]
We get 
 \[
\ \begin{aligned}
\ & | \mathcal  C^1_\psi(H_{i,m} dA) (w) - \mathcal  C^1_\psi\nu_i (w) | \\
\  \le & \left | \dfrac{m^2}{\pi} \int _{|\lambda-w|\ge \frac{1}{\sqrt{m}}} \int _{|z-\lambda | < \frac{1}{m}} \dfrac{\psi (z) - \psi (\lambda )}{z - w} dA(z) d\nu_i (\lambda )\right |+ 2 C_0|\nu_i | ( B(w, \frac{1}{\sqrt{m}}))\\
\  \le &  \dfrac{m^2}{\pi} \int _{|\lambda-w|\ge \frac{1}{\sqrt{m}}} \int _{|z-\lambda | < \frac{1}{m}} \dfrac{C_0|z - \lambda |}{|w - \lambda | - |z - \lambda |} dA(z) d\nu_i (\lambda )+ 2 C_0|\nu_i | ( B(w, \frac{1}{\sqrt{m}}))\\
\  \le & C_0\dfrac{\frac{1}{m}}{\frac{1}{\sqrt{m}} - \frac{1}{m}}|\nu_i | ( B(w, \frac{1}{\sqrt{m}})^c) + 2 C_0|\nu_i | ( B(w, \frac{1}{\sqrt{m}}))\\
\le & \frac{C_0}{\sqrt{m} - 1} \|\nu_i \| + 2C_0 |\nu_i | ( B(w, \frac{1}{\sqrt m})).
\end{aligned}
\]
Therefore,
 \[
 \ \lim_{m\rightarrow \infty} \mathcal  C^1_\psi(H_{i,m} dA) (w) = \mathcal  C^1_\psi\nu_i (w)
 \]
for $w\notin A.$ For $\lambda _0 \in \Omega^G_N$ and $\epsilon > 0,$ we can choose a $\delta > 0$ and $m_0$ such that 
\[
\begin{aligned}
\ & | \mathcal  C^1_\psi(H_{i,m} dA) (w) - \mathcal  C^1_\psi\nu_i (w) |  \\
\le & 2C_0 |\nu_i | ( B(w, \frac{1}{\sqrt m})) + \frac{C_0}{\sqrt{m} - 1} \|\nu_i \| \\
\le & 2C_0 |\nu_i | ( B(\lambda_0, \delta + \frac{1}{\sqrt m})) + \frac{C_0}{\sqrt{m} - 1} \|\nu_i \| \\
< &\epsilon 
\end{aligned}
\]
where $w\in B(\lambda _0, \delta )\setminus A$ and $m\ge m_0.$ Since $\mathcal  C^1_\psi\nu_i (w)$ is continuous at $\lambda _0,$ $\delta$ can be chosen to ensure
 \[
 \ | \mathcal  C^1_\psi\nu_i (w) - \mathcal  C^1_\psi\nu_i (\lambda _0) | < \epsilon
 \]
where $w\in B(\lambda _0, \delta )\setminus A.$ 
It is easy to verify that $\mathcal  C^1_\psi(H_{i,m} dA)$ is a smooth function. For $k > 1,$ clearly $\mathcal  C^k_\psi\nu_i (w)$ is a smooth function. Define
  \[
 \ M_N^{Gm}(z) = \left[ \begin{array}{cccc}
\mathcal  C^1_\psi(H_{1,m} dA), & \mathcal  C^1_\psi(H_{2,m} dA), &..., & \mathcal  C^1_\psi(H_{N-1,m} dA)\\
\mathcal  C^2_\psi(\nu_1), & \mathcal  C^2_\psi(\nu_2), &..., & \mathcal  C^2_\psi(\nu_{N-1})\\
..., & ..., &..., & ... \\
\mathcal  C^{N-1}_\psi(\nu_1), & \mathcal  C^{N-1}_\psi(\nu_2), &..., & \mathcal  C^{N-1}_\psi(\nu_{N-1}) \end{array} \right].
 \]
We can choose $\epsilon$ small enough so that 
 \[
 \ M_N^{Gm}(w), ~ M_N^G(w)
 \]
are invertible for $w\in B(\lambda _0, \delta )\setminus A$ and $m > m_0.$ 
Define
 \[
 \ B_N^{Gm}(z) = \left[ \begin{array}{c}
\mathcal  C^1_\psi(H_{N,m} dA)\\
\mathcal  C^2_\psi(\nu_{N})\\
... \\
\mathcal  C^{N-1}_\psi(\nu_{N}) \end{array} \right],
 \]
 \[
 \ A_N^{Gm}(z) = [ \mathcal C(H_{1,m} dA), \mathcal C(H_{2,m} dA), ... , \mathcal C(H_{N - 1,m} dA) ]
 \]
and
 \[
 \ c_N^{Gm}(z) = \mathcal C(H_{N,m} dA)(z).
 \]

For a smooth function $\phi$ with compact support in $B(\lambda _0, \delta ),$ using the definition (3-8) and Lebesgue's Dominated Convergence Theorem, we get the following calculation,
\[
\begin{aligned}
\ & \int \bar\partial \phi (z) a_N(z) dA(z) \\
\ = & \lim_{m\rightarrow \infty} \int \bar\partial \phi (z) \left ( ( M_N^{Gm}(z) )^{-1} B_N^{Gm} (z) \right ) dA(z) \\
\ =& - \lim_{m\rightarrow \infty} \int \phi (z) \bar\partial \left ( ( M_N^{Gm}(z) )^{-1} B_N^{Gm} (z) \right ) dA(z) \\
\ =& \lim_{m\rightarrow \infty} \int \phi (z) ( M_N^{Gm}(z) )^{-1} \left ((\bar\partial M_N^{Gm}(z)) ( M_N^{Gm}(z) )^{-1} B_N^{Gm} (z)  - \bar\partial B_N^{Gm} (z)\right ) dA(z).
\end{aligned} \tag{3-12}
\]
On the other hand,
\[
 \ \bar\partial M_N^{Gm}(z) = \bar\partial\psi(z)\left[ \begin{array}{cccc}
- \mathcal C(H_{1,m} dA), & - \mathcal C(H_{2,m} dA), &..., & - \mathcal C(H_{N-1,m} dA)\\
-2 \mathcal  C^1_\psi(\nu_1), & -2 \mathcal  C^1_\psi(\nu_2), &..., & -2 \mathcal  C^1_\psi(\nu_{N-1})\\
..., & ..., &..., & ... \\
- (N-1) \mathcal  C^{N-2}(\nu_1), & - (N-1) \mathcal C^{N-2}(\nu_2), &..., & - (N-1) \mathcal C^{N-2}(\nu_{N-1}) \end{array} \right].
 \]
Therefore,
\[
 \ (\bar\partial M_N^{Gm}(z)) ( M_N^{Gm}(z) )^{-1} = - \bar\partial\psi(z) \left[ \begin{array}{ccccc}
 A_N^{Gm}(z) & &( M_N^{Gm}(z) )^{-1}, & &\\
2 , & 0, &..., & 0, & 0 \\
..., & ..., &..., & ..., & ... \\
0 , & 0, &..., & N - 1, & 0 \end{array} \right].
 \]
Hence,
 \[
 \ (\bar\partial M_N^{Gm}(z)) ( M_N^{Gm}(z) )^{-1} B_N^{Gm} (z) - \bar \partial B_N^{Gm} (z) = - \bar\partial\psi(z)\left[ \begin{array}{c}
 A_N^{Gm}(z) ( M_N^{Gm}(z) )^{-1} B_N^{Gm} - c_N^{Gm} \\
 0 \\
... \\
0  \end{array} \right].
 \]
Using (3-11), we see that
 \[
 \ \lim_{m\rightarrow \infty} \left ( A_N^{Gm}(z) ( M_N^{Gm}(z) )^{-1} B_N^{Gm} - c_N^{Gm} \right ) = 0 ~ a.e. ~Area |_{B(\lambda _0, \delta )}.
 \]
Since each component of the above vector function is less than
 \[
 \ M \int \dfrac{1}{|w - z | } d |\nu_i| (z) ~a.e.~ Area |_{B(\lambda _0, \delta )}, 
 \]
applying Lebesgue's Dominated Convergence Theorem to the last step of (3-12), we conclude
 \[
 \ \int \bar\partial \phi (z) a_N(z) dA(z) = 0.
 \]
By Weyl's lemma, we see that $a_N(z)$ is analytic on $B(\lambda _0, \delta ).$ From equation (3-8), we get
\[
 \ \mathcal  C^1_\psi \langle F_N, G\rangle \mu_1)(z) = \sum_{k=1}^{N-1} a_{kj}(z)\mathcal  C^1_\psi \langle F_k, G\rangle \mu_1)(z), ~ a.e.~ Area |_{ B(\lambda _0, \delta )}.
 \]
The above equation implies (3-9) since 
 \[
 \ \bar \partial \mathcal  C^1_\psi(\nu_i)(z) = - \mathcal C(\nu_i)(z) ~ a.e. ~ Area. 
 \]
For equation (3-10), let $\phi$ be a smooth function with compact support in $B(\lambda _0, \delta )$ and let $\nu$ be a compactly supported finite measure, we get 
\[
 \  \int \bar \partial \phi (z) \mathcal C\nu (z) dA(z) = \pi \int \phi (z) d\nu (z).
\]
Apply the above equation to the both sides of the equation (3-9) for $j=N$ and using
 \[
 \ \bar \partial \phi (z) a_{kj} (z) = \bar \partial (\phi (z) a_{kj} (z)), ~ z\in B(\lambda _0, \delta ),
 \]
we conclude
 \[
 \ \int \phi \langle F_N, G \rangle d\mu_1 = \int \phi \sum_{k=1}^{N-1} a_{kj}\langle F_k, G \rangle d\mu_1 .
 \]
Hence the equation (3-10) follows. This completes the proof of the lemma.

\begin{Corollary} \label{N1ABPECorollary}
Let $G,$ $\Omega^G,$ and $\Omega^G_i$ be as in Lemma \ref{KeyLemma}. Suppose $G\perp\mathcal K^\psi_{N-1}$ (satisfies (3-4)). Then $\Omega^G_i\subset U_{N-1}\cup U_N.$
\end{Corollary}

{\bf Proof:} Without loss of generality, we assume that $j = N.$ From Lemma \ref{KeyLemma}, for $\lambda _0\in \Omega^G_N,$ there exists $\delta >0$ such that $B(\lambda _0,\delta)\subset Int(\sigma(S))$ and the equations (3-9) and (3-10) hold, which imply (3-3). For $r_1,r_2,...,r_N\in Rat(\sigma(S)),$ let
 \[
 \ F = \sum_{i=1}^N r_iF_i.
 \]
Notice that
 \[
 \ r_i(\lambda)\mathcal C_\psi^k\langle F_i, G \rangle\mu_1) = \mathcal C_\psi^k\langle r_iF_i, G \rangle\mu_1)
 \]
since $G\perp\mathcal K^\psi_{N-1}.$
Then
 \[
 \ \sum_{i=1}^N r_i(\lambda)\mathcal C_\psi^k(\langle F_i, G \rangle\mu_1) (\lambda) = \mathcal C_\psi^k(\langle F, G \rangle\mu_1) (\lambda),
 \] 
for $k = 1,2,...,N-1.$ Now using the equation (3-9) for $\lambda\in B(\lambda _0, \delta )\setminus A,$ we get
\[
 \ \sum_{i=1}^{N-1} (r_i(\lambda) + a_{Ni} (\lambda ) r_N(\lambda)) \mathcal C_\psi^k(\langle F_i, G \rangle\mu_1) (\lambda) = \mathcal C_\psi^k(\langle F, G \rangle\mu_1) (\lambda),
 \]
equivalently,
 \[
 \ M_N^G(\lambda ) \left[ \begin{array}{c}
 r_1(\lambda) + a_{N1} (\lambda ) r_N(\lambda) \\
 r_2(\lambda) + a_{N2} (\lambda ) r_N(\lambda) \\
...,  \\
 r_{N-1}(\lambda) + a_{N,N-1} (\lambda ) r_N(\lambda) \end{array} \right] = \left[ \begin{array}{c}
 \mathcal C_\psi^1(\langle F, G \rangle\mu_1) (\lambda) \\
 \mathcal C_\psi^2(\langle F, G \rangle\mu_1) (\lambda)  \\
...,  \\
 \mathcal C_\psi^{N-1}(\langle F, G \rangle\mu_1) (\lambda) \end{array} \right].
 \]
where the inverse of $M_N^G(\lambda )$ is bounded on $B(\lambda _0, \delta )\setminus A$ and $a_{Ni}$ are analytic on $B(\lambda _0, \delta ).$
Therefore, there exists a positive constant $M$ such that
 \[
 \ \sup_{1\le k \le N-1, \lambda \in B(\lambda _0, \frac{\delta}{2} )}|r_{k}(\lambda) + a_{Nk} (\lambda ) r_N(\lambda) | \le M \|F\|,
\]
which implies (3-2).
Hence, Lemma 3.1 implies $\Omega_N^G\subset U_{N-1}\cup U_N.$ 

Now let us recursively construct other sets such as $\Omega_{ij}^G$ for a given $G\perp\mathcal K^\psi_{N-1}.$ We will only describe the algorithm for $k = N - 2$ and the other cases will follow recursively. Let $E_N^G = \Omega^G$ and $E_{N-1}^G = \cup _{i=1}^N \Omega_i^G.$ Let $M_{ij}^G$ be an $N-2$ by $N-2$ submatrix of $M^G$ by eliminating the first two rows and the $i$ and $j$ columns. Define
 \[
 \ \Omega_{ij}^G = (Int(\sigma(S))\cap A^c \cap \{ z: | det(M_{ij}^G(z)) | > 0 \}) \setminus clos(E_N^G\cup E_{N-1}^G).
 \]
Without loss of generality, let us assume that $i = N-1$ and $j = N.$
Similar to Lemma \ref{KeyLemma}, one can prove that for $\lambda _0\in \Omega_{N-1,N}^G,$ there exist $\delta >0,$ analytic functions $a_i(z)$ and $b_i(z)$ on $B(\lambda _0, \delta )\subset Int(\sigma(S))$ such that
 \[
 \ F_{N-1} = \sum _{i=1}^{N-2} a_i(z) F_i(z), ~ F_N = \sum _{i=1}^{N-2} b_i(z) F_i(z), ~ a.e.\mu_1 | _{B(\lambda _0, \delta )},\tag{3-13}
 \]
and there exists a constant $M > 0$ such that
 \[
 \ \sup_{1\le k \le N-2, \lambda \in B(\lambda _0, \frac{\delta}{2} )}|r_{k}(\lambda) + a_k (\lambda ) r_{N-1}(\lambda) + b_k (\lambda ) r_N(\lambda) | \le M \| F\|,\tag{3-14}
\]
where $r_1,r_2,...,r_N\in Rat(\sigma(S))$ and $F = \sum_{i=1}^N r_iF_i.$ (3-13) and (3-14) are the same as (3-2) and (3-3) for the case $k=N-2.$ Let
 \[
 \ E_{N-2}^G = \cup _{i < j}^N \Omega_{ij}^G.\tag{3-15}
 \]

\begin{Corollary}\label{N2ABPECorollary}
Let $E_{N-2}^G$ be as in (3-15). Suppose $G\perp\mathcal K^\psi_{N-1}$ (satisfies (3-4)). Then
 \[
 \ E_{N-2}^G \subset U_{N-2} \cup U_{N-1} \cup U_N.
 \]
\end{Corollary}

The proof is the same as Corollary \ref{N1ABPECorollary}. Therefore we can recursively construct $E_k^G$ for $k = 1,2,...,N$ such that 
 \[
 \ E_k^G \subset \bigcup_{i=k}^N U_i\tag{3-16}
 \]
where the proof for $k=N$ is from Lemma \ref{KeyLemma0}, $k=N-1$ is from Corollary \ref{N1ABPECorollary}, and $k=N-2$ is from Corollary \ref{N2ABPECorollary}.

The following theorem proves, under the conditions $S$ satisfies the property $(N,\psi ),$ the set $\cup_{k=1}^N E_k^G$ is big.

\begin{Theorem} \label{BigETheorem}
Let $E_i^G$ be constructed for $i = 1,2,...,N$ as above. Suppose $\{G_j\}\subset (\mathcal K^\psi_{N-1})^\perp$ is a dense subset, then 
 \[
 \ spt \mu _1 \subset clos\left (\bigcup _{i=1}^N \bigcup _{j=1}^\infty E_i^{G_j}\right).
 \]
\end{Theorem} 

{\bf Proof:} First we prove 
 \[
 \ \mu_1 \left ( Int(\sigma (S)) \setminus clos\left (\bigcup _{i=1}^N \bigcup _{j=1}^\infty E_i^{G_j}\right)\right ) = 0.
 \]
Suppose that $B(\lambda _0, \delta) \subset Int(\sigma (S))$ and $B(\lambda _0, \delta)\cap clos\left (\bigcup _{i=1}^N \bigcup _{j=1}^\infty E_i^{G_j}\right) = \emptyset,$ then by construction of $E_i^{G_j},$ we conclude that 
 \[
 \ \mathcal C_\psi^{N-1}(\langle F_i,G_j\rangle\mu_1) (z) = 0
 \]
on $B(\lambda _0, \delta ),$ where $i = 1,2,...,N.$ By taking $\bar \partial$ in the sense of distribution, we see that
 \[
 \ \mathcal C(\langle F_i,G_j\rangle\mu_1) (z) = 0
 \]
a.e. Area on $B(\lambda _0, \delta )$ since $Area(\{\bar\partial\psi = 0\}\cap \sigma (S)) = 0,$ where $i = 1,2,...,N.$ 
For a smooth function $\phi$ with compact support in $B(\lambda _0, \delta ),$
 \[
 \ \int \phi (z) \langle F_i,G_j\rangle d\mu_1 = \dfrac{1}{\pi}\int \bar \partial \phi(z) \mathcal C(\langle F_i,G_j\rangle\mu_1) (z) dA(z) =0.
\]
Therefore, 
 \[  
 \ \langle F_i(z),G_j(z)\rangle = 0. ~ a.e. ~ \mu_1 | _{B(\lambda _0, \delta)} \tag{3-17}
 \]
where $i = 1,2,...,N.$ From (1-4), we see that for $P\in \oplus _{k=1}^m L^2(\mu _k |_{B(\lambda _0, \delta)}),$ (3-17) implies $(P,G_j) = 0.$ Therefore, 
 \[
 \ \oplus _{k=1}^m L^2(\mu _k |_{B(\lambda _0, \delta)})\subset \mathcal K^\psi_{N-1}.
 \]
Hence, $\mu_1 | _{B(\lambda _0, \delta)} = 0$ since $M_z|_{\mathcal K^\psi_{N-1}}$ is pure.

Now assume $B(\lambda _0, \delta)\cap clos(Int(\sigma(S))) = \emptyset.$ For $N>1,$ the function $\mathcal C_\psi^{N-1}(\langle F_i,G_j\rangle\mu_1) (z)$ is continuous on $\mathbb C \setminus A$ and is zero on $\mathbb C \setminus \sigma(S).$ Hence,
\[
 \ \mathcal C_\psi^{N-1}(\langle F_i,G_j\rangle\mu_1) (z) = 0
 \]
on $B(\lambda _0, \delta )\setminus A,$ where $i = 1,2,...,N.$ Using the same proof as above, we see that $\mu_1 | _{B(\lambda _0, \delta)} = 0.$ This implies $spt\mu_1\subset clos(Int(\sigma (S))).$ The theorem is proved.

{\bf Proof of Theorem \ref{MTheorem2}:} From (3-16) and Theorem \ref{BigETheorem}, we get
 \[
 \ \bigcup_{i = 1}^N \partial U_i \subset \sigma _ e (S) \subset spt(\mu_1) \subset clos\left(\bigcup_{i = 1}^N U_i\right).
 \] 
This implies
 \[
 \ \sigma _ e (S) = \bigcup_{i = 1}^N \partial U_i
 \] 
since $\sigma _ e (S)\cap U_i =\emptyset.$ This completes the proof.

For a positive finite measure $\mu$ with compact support on $\mathbb C, $ definite
 \[
 \ P^2(\mu | 1,\bar z,...,\bar z^{N-1}) = clos\{p_1(z) + p_2(z) \bar z +...+ p_N(z) \bar z^{N-1}:~ p_1,p_2,..., p_N\in \mathcal P\}
 \]
and $S_{N,\mu}$ as the multiplication by $z$ on $P^2(\mu | 1,\bar z,...,\bar z^{N-1}).$ Then $S_{N,\mu}$ is a multicyclic subnormal operator with the minimal normal extension $M_\mu,$ the multiplication by $z,$ on $L^2(\mu ).$

\begin{Corollary} \label{SimpleCorollary}
Suppose that $S_{2,\mu}$ on $P^2(\mu | 1, \bar z, \bar z^2)$ is pure, then the operator $S_{1,\mu}$ on $P^2(\mu | 1, \bar z)$ satisfies
 \[
 \ \sigma (S_{1,\mu})  = clos(\sigma (S_{1,\mu}) \setminus \sigma _e(S_{1,\mu})).
 \]
\end{Corollary}

Proof: Since   
 \[
 \ \mathcal K^{\bar z}_1 = clos (span (\bar z^k P^2(\mu |1,\bar z): 0 \le k \le 1)) = P^2(\mu | 1, \bar z, \bar z^2)
 \]
and $S_{2,\mu}$ on  $P^2(\mu | 1, \bar z, \bar z^2)$ is pure. Therefore, the result follows from Theorem \ref{MTheorem2}.

It seems strong to assume that $S_{2,\mu}$ on $P^2(\mu | 1, \bar z, \bar z^2)$ is pure in the corollary. We believe that the condition can be reduced to assume that $S_{1,\mu}$ on $P^2(\mu | 1, \bar z)$ is pure. However, we are not able to prove the result under the weaker conditions. We will leave it as an open problem for further research. 

\begin{Problem} 
Does Corollary \ref{SimpleCorollary} hold under the weaker assumption that $S_{1,\mu}$ on $P^2(\mu | 1, \bar z)$ is pure?
\end{Problem}

\begin{Corollary} \label{SpecialCorollary}
Let $S$ on $\mathcal H$ be a pure rationally $N-$cyclic subnormal operator with $\mathcal H = R^2(S|F_1,F_2,...,F_N)$ and let $M_z$ be its minimal normal extension on $\mathcal K$ satisfying (1-1) to (1-4). Suppose that there exists a smooth function $\psi$ on $\mathbb C$ such that $Area(\{\bar \partial \psi = 0\}\cap \sigma(S) ) = 0$ and $\psi(M_z)\mathcal H\subset \mathcal H.$ Then there exist bounded open subsets $U_i$ for $ 1\le i \le N$ such that 
\[
 \ \sigma _e (S) = \bigcup _{i=1}^N \partial U_i, ~\sigma (S) \setminus \sigma _e(S) = \bigcup _{i=1}^N U_i,
 \]
and 
 \[
 \ dim ker (S - \lambda )^* = i.
 \]
for $\lambda \in U_i.$
\end{Corollary}

Notice that Example \ref{FRSCExample} and \ref{TCExample} are special cases of Corollary \ref{SpecialCorollary}. It seems that further results could be obtained for the special cases where $S$ satisfies the conditions of Corollary \ref{SpecialCorollary}. Moreover, we might be able to combine the methodology in \cite{my97} to obtain the structural models for the class of subnormal operators, which might extend  Xia's model for subnormal operators with finite rank self-commutators.

\begin{Problem}
Can the structure of subnormal operators in Corollary \ref{SpecialCorollary} be characterized?
\end{Problem}

\bibliography{Bibliography}
\end{document}